\documentclass{article}

\usepackage{arxiv}

\usepackage[utf8]{inputenc} 
\usepackage[T1]{fontenc}    
\usepackage{hyperref}       
\usepackage{url}            
\usepackage{booktabs}       
\usepackage{amsfonts}       
\usepackage{nicefrac}       
\usepackage{microtype}      
\usepackage{lipsum}		
\usepackage{graphicx}
\usepackage{natbib}
\usepackage{doi}

\usepackage{enumitem}		
\usepackage{indentfirst}		
\usepackage{color}			
\usepackage{graphicx}			
\usepackage{microtype} 
\usepackage{scalerel}[2016/12/29]
\usepackage{mathtools,amsfonts,amssymb,amsxtra,empheq,amsthm}
\usepackage{comment}
\usepackage{stackengine}

\newcommand{\disp}{\displaystyle}
\newcommand{\non}{\nonumber}

\newcommand{\defi}{\coloneqq}

\makeatletter
\newcommand*{\bigcdot}{}
\DeclareRobustCommand*{\bigcdot}{%
  \mathbin{\mathpalette\bigcdot@{}}%
}
\newcommand*{\bigcdot@scalefactor}{.5}
\newcommand*{\bigcdot@widthfactor}{1.5}
\newcommand*{\bigcdot@}[2]{%
  \sbox0{$#1\vcenter{}$}
  \sbox2{$#1\cdot\m@th$}%
  \hbox to \bigcdot@widthfactor\wd2{%
    \hfil
    \raise\ht0\hbox{%
      \scalebox{\bigcdot@scalefactor}{%
        \lower\ht0\hbox{$#1\bullet\m@th$}%
      }%
    }%
    \hfil
  }%
}
\makeatother

\newcommand{\N}{\mathbb{N}}
\newcommand{\Z}{\mathbb{Z}}

\newcommand{\R}{\mathbb{R}}

\newcommand{\T}{\mathbb{T}}

\newcommand{\Ttil}{\widetilde{\mathbb{T}}}

\newcommand{\deri}{\vartriangle}

\newcommand{\bigsigma}{\text{S}(\T)}
\newcommand{\bigsigmatil}{\text{S}(\widetilde{\T})}
\newcommand{\bigsigman}{\text{S}(\mathbb{T}_n)}
\newcommand{\sigmadelta}{\text{S}_{P_\delta}(\T)}

\newcommand{\J}{I(\mathbb{T})}

\newtheorem{teorema}{Theorem}[section]
\newtheorem{proposicao}[teorema]{Proposition}
\newtheorem{definicao}[teorema]{Definition}
\newtheorem{lema}[teorema]{Lemma}
\newtheorem{corolario}[teorema]{Corollary}
\newtheorem{exemplo}{Example}

{\theoremstyle{definition}
\newtheorem{afirmacao}{Affirmation}[teorema]

\newtheorem*{prova}{Proof}}

\title{A Different Demonstration for Integral Identity Across Distinct Time Scales}

\author{ \href{https://orcid.org/0000-0002-6800-3289}{\includegraphics[scale=0.06]{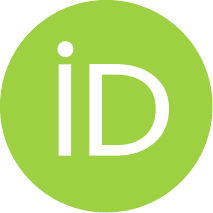}\hspace{1mm}Patrick Oliveira} \\
	Departamento de Matemática\\
	Universidade Federal de Minas Gerais\\
	\texttt{patrickoliveira@mat.dout.ufmg.br} \\
}

\renewcommand{\shorttitle}{\textit{arXiv} Template}

\hypersetup{
pdftitle={A Carathéodory Extension Theorem for Families Simpler Than an Algebras},
pdfsubject={q-mat.NC},
pdfauthor={Patrick Oliveira},
pdfkeywords={Time Scales, Integral Identity},
}

\begin{document}
\maketitle

\begin{abstract}

In the theory of time scales, given $\T$ a time scale with at least two distinct elements, an integration theory is developed using ideas already well known as Riemann sums. Another, more daring, approach is to treat an integration theory on this scale from the point of view of the Lebesgue integral, which generalizes the previous perspective. A great tool obtained when studying the integral of a scale $\T$ as a Lebesgue integral is the possibility of converting the ``$\Delta$-integral of $\T$'' to a classical integral of $\R$. In this way, we are able to migrate from a calculation that is sometimes not so intuitive to a more friendly calculation. A question that arises, then, is whether the same result can be obtained just using the ideas of integration via Riemann sums, without the need to develop the Lebesgue integral for $\T$. And, in this article, we answer this question affirmatively: In fact, for integrable functions an analogous result is valid by converting a $\Delta$-integral over $\T$ to a riemannian integral of $\R$.
\end{abstract}


\section{Introduction}

A famous and very useful result within the theory of time scales is Theorem 5.2. of \cite{expression} which guarantees that, given $\T$ a time scale, $E \subset \T$ a $\Delta_\T$-measurable set and $f:\T \longrightarrow \R$ a $\Delta_\T$-integrable function, then the $\Delta_\T$-Lebesgue integral of $f$ over $E$, denoted by $\int_{E} f(s) \, \Delta_\T s $, satisfies the identity
\begin{eqnarray}
\int_{E} f(s) \, \Delta_\T s & = & \int_{E} f(s) \, ds + \sum_{\tau < \sigma_\T(\tau)}^{\tau \in E} f(\tau) \cdot \mu_\T(\tau),
\label{eq_1}
\end{eqnarray}

\noindent where $\int_{E} f(s) \, ds$ denotes the Lebesgue integral over $\R$ (\cite{bartle, halmos}) and the maps $\sigma_\T$ and $\mu_\T$ represent, respectively, the foward jump operator and the graininess of the scale $\T$ (\citep{dynamic, advances}). Our investigation, in this article, is whether this same result can be demonstrated for the Riemanian $\Delta_\T$-integrals, without the need to use the entire Lebesgue integral tool. And, in fact, using Abel's Sum Lemma (Lemma \ref{lema_abel}), we are able to demonstrate that a formula analogous to equation (\ref{eq_1}) holds for $\Delta_\T$-integrable functions (Corollary \ref{cor_fim}). From this result, we developed some examples, corollaries and, in particular, an application on the convergence of sequences of integrals when the integrand is fixed and the time scale varies.

\section{Basic Results and Definitions}
\label{sec:1}

Let us establish some common notations and hypotheses throughout the article. We will denote by $\T$ and $\Ttil$ time scales such that $\T \subset \Ttil$. We agree that $\#\T \geq 3$, that is, $\T$ has at least three distinct points, hence the construction of the $\Delta_\T$-integral over $\T$ is well defined. 

Given functions $f: \Ttil \longrightarrow \R$ and $g: \T \longrightarrow \R$, respectively, $\Delta_{\Ttil}$ and $\Delta_\T$-integrable in the Riemannian sense (see \cite{advances}, page 117), we represent, respectively, the $\Delta_{\Ttil}$-integral of $f$ over $[a,b]\cap\Ttil$ and the $\Delta_\T$-integral of $g $ over $[a,b] \cap \T$, by
\begin{equation}
\int_{a}^{b} f(s) \, \Delta_{\Ttil} s \quad \mbox{and} \quad \int_{a}^{b} g(s) \, \Delta_\T s,
\end{equation}

\noindent where $a,b \in \T$, with $a < \rho_\T(b)$, with $\rho_\T$ being the backward jump operator of $\T$ (\cite{dynamic}).


Below we make use of the main definitions and constructive results of this article.

\begin{definicao}

Let $\delta > 0$ and $a,b \in \T$. For each delta-partition $P_\delta = \{a = t_1, ... , t_n = b\} \in \mathcal{P}_\delta(a,b)$ of the interval $[a,b]_ \T \defi [a,b] \cap \T$, we define the $P_\delta$-partition function $\sigmadelta:[a,b] \longrightarrow \mathbb{R}$ by
\begin{equation} \label{Tdelta_first_definition}
\sigmadelta(t) \; \defi \;  \left\{\begin{array}{rl} 
\hspace{0.25cm} t_i, & \hspace{0.25cm} \emph{if} \hspace{0.25cm}  t \in [t_{i-1},t_i) \hspace{0.25cm} \emph{and} \hspace{0.25cm} i \in \{2, \cdots, n-1\} \\ & \\
b, & \hspace{0.25cm} \emph{if} \hspace{0.25cm} t \in [t_{n-1},b],\\
\end{array}\right.
\end{equation}


\noindent or, equivalently, $\sigmadelta(t) = \inf\{ s \in P_\delta : s > t\}$ for all $t \in [a,b)$ and $\sigmadelta(b) = b$ .

\label{def_sigmadelta}
\end{definicao}

\begin{definicao}

Let $a,b \in \T$. We define the set $\J \subset [a,\rho_\T(b))$ as the meeting of all real intervals $[\tau, \sigma_\T(\tau))$ such that $\tau \in \T$ is right discrete in $\T$ and $a \leq \tau < \rho_\T(b)$. That way,
\begin{equation*}
\J \; \defi \; \hspace{-0.15cm} \displaystyle\stackrel{\tau \hspace{0.1cm} \in \hspace{0.1cm} [a,\rho_\T(b))_\mathbb{T}}{\bigcup_{\sigma_\T(\tau)>\tau}} \hspace{-0.15cm} [\tau,\sigma_\T(\tau)).
\end{equation*}

\label{def_J}
\end{definicao}


\vspace{-0.15cm}

Note that $\J \subset [a,\rho_\T(b))$. Indeed, if $t \in \J$, then there exists right-discrete point $\tau \in [a,\rho_\T(b))_\T$ in $\T$ such that $t \in [\tau, \sigma_\T(\tau))$. Therefore, $a \leq \tau \leq t < \sigma_\T(\tau)$. As $\tau < \rho_\T(b)$, then $\sigma_\T(\tau) \leq \rho_\T(b)$. Therefore $a \leq t < \sigma_\T(\tau) \leq \rho_\T(b)$, where $t \in [a,\rho_\T(b))$. Note that $t \in \J $ if and only if there exists $\tau \in [a,\rho_\T(b))_\T$ right-discrete in $\T$ such that $t \in [\tau, \sigma_\T(\tau)) \subset \J$. Therefore, for the sake of simplicity, sometimes when we say that $t \in \J$, we will abbreviate the previous equivalence by saying $t \in [\tau,\sigma_\T(\tau)) \subset \J$ . It is agreed that $\tau \in \T$ is a discrete point to the right of $\T$.

\vspace{0.15cm}

\begin{exemplo}

For suitable $a,b \in \R$, consider the following scales

\begin{enumerate}
\item[]\textbf{\emph{(i)}} 
If $\T = \R$, then $I(\R) = \emptyset$, since every point of $\R$ is right-dense.

\item[]\textbf{\emph{(ii)}} 
If $\T = \Z$, taking $a,b \in \Z$ such that $a < b-1$, then $I(\Z) = [a,a+1) \cup [a+1 ,a+2) \cup \cdots \cup [b-2,b-1)$, that is, $I(\Z) = [a,b-1)$.

\item[]\textbf{\emph{(iii)}} 
If $\T = \{0\} \cup \{q^n : n \in \N\}= \overline{q^{\N}}$, where $q \in (0,1)$ , taking $a=0$ and $b = q$ we have \\ $I(\overline{q^{\N}}) = [q^3,q^2)\cup[q^4,q^3)\cup \cdots \cup [q^{n+1},q^n) \cup \cdots$, then, $I(\overline{q^{\N}}) = (0,q^2)$.
\end{enumerate}

\end{exemplo}

\vspace{0.15cm}

\begin{proposicao}

Let $\{\sigmadelta\}_{\delta > 0}$ be the family formed by all functions $\sigmadelta$
defined according to \emph{Definition \ref{def_sigmadelta}} for $a,b \in \T$, with $a < \rho_\T(b)$, and all $\delta > 0$. Then $\displaystyle\lim_{\delta \to 0} \sigmadelta = \bigsigma$, where $\bigsigma: [a,b] \longrightarrow \mathbb{R}$ is defined by

\begin{equation} \label{bigsigma_definition}
\bigsigma(t) \; \defi \; \left\{\begin{array}{rl}
\sigma_\T(\tau), & \emph{if} \hspace{0.25cm} t \in [\tau,\sigma_\T(\tau)) \subset \J   \\ & \\
t, & \emph{if} \hspace{0.25cm} t \in [a,\rho_\T(b))\setminus \J \\ & \\
b, & \emph{if} \hspace{0.25cm} t \in [\rho_\T(b), b]
\end{array}\right.
\end{equation}


\noindent for all $t \in [a,b]$, where $\J \subset [a,\rho(b))$ is defined for $[a,b]_\T$ according to \emph{Definition \ref{def_J}} (equivalently, we have $\bigsigma(t) = \inf\{s \in [a,b]_\T : s > t\}$, if $t \in [a,\rho_ \T(b))$, and $\bigsigma(t) = b$, if $t \in [\rho_\T(b), b]$). Furthermore, the family of functions $\sigmadelta$ converges uniformly to $\bigsigma$ when $\delta$ tends to $0^+$ in the following sense:

\hspace{3.75cm}\begin{minipage}[h]{10cm}

Given $\varepsilon > 0$, there exists $\delta_0 > 0$ such that, for all $P_\delta \in \mathcal{P}_\delta(a,b)$ and all $t \in [a, b]$, we have $|\sigmadelta(t) - \bigsigma(t)| < \varepsilon$ for all $\delta \in (0,\delta_0)$.
\end{minipage}

\label{prop_bigsigma}
\end{proposicao}

\newpage

\begin{prova}

Throughout the proof we will denote, for simplicity of notation, $\rho = \rho_\T$ and $\sigma = \sigma_\T$. We will show that there is $\delta_0 > 0$ such that for all $\delta \in (0, \delta_0)$ we have $|\sigmadelta(t) - \bigsigma(t)| < \delta$ for every $t \in [a,b]$ and for every $P_\delta$ delta-partition of $[a,b]_\T$.


Given $t \in [a,b]$ and a partition $P_\delta$, to limit the module $|\sigmadelta(t)-\bigsigma(t)|$, we need to determine the values of the functions $\sigmadelta$ and $\bigsigma$ in $t \in [a,b]$. Therefore, we will use the following statement:

\vspace{0.5cm}

\begin{afirmacao}

There exists $\delta_0 > 0$ such that for every $\delta \in (0,\delta_0)$, every delta-partition $P_\delta \in \mathcal{P}(a,b)$ and every $t \in [a,\rho(b))$, there are $t_{i-1},t_i \in P_\delta$ such that $t \in [t_{i-1},t_i)\subset [a,\rho(b))$, with $t_{i} \leq \rho(b)$.
\end{afirmacao}

\begin{prova}

If $b \in \T$ is left discrete in $\T$, let $\delta_0 > 0$ given by $\delta_0 = (b - \rho(b))/2$. Note that in this case $\delta_0 = \mu(\rho(b))/2$. Thus, every partition $P_\delta$ with $\delta \in (0, \delta_0)$ will be in the form $P_\delta = \{a = t_1, t_2, ..., t_{n-2}, t_{n-1} = \rho(b), t_n = b\}$. Thus, if $t \in [a,\rho(b))$, then $t \in [t_{i-1},t_i)$ for some $t_{i-1}, t_i \in P_\delta$, with $t_i$ being at most $t_{n-1} = \rho(b)$.


If $b \in \T$ is left dense in $\T$, let $\delta_0 = 1$. In this case, for every delta-partition $P_\delta = \{a = t_1, ..., t_{n-1}, t_n =b\}$, we have that, if $t \in [a,\rho (b)) = [a,b)$, then there are $t_{i-1},t_i \in P_\delta$ such that $t \in [t_{i-1}, t_i)$. Clearly $t_i$ can be at most $t_n = b = \rho(b)$. \hfill $\vartriangle$
\end{prova}


We will use this statement in items (i) and (ii) below. Choosing $\delta_0 > 0$ suitable, if
$\delta \in (0,\delta_0)$. Let's show for every delta-partition $P_\delta$ we have
\begin{equation*}
|\sigmadelta(t) - \bigsigma(t)| =  \left\{\begin{array}{ll} 
|t_i - \sigma(\tau)|, & \emph{if} \hspace{0.25cm} t \in [\tau,\sigma(\tau)) \subset  \J \subset [a,\rho(b))\\ & \\ 
|t_i - t|,            & \emph{if} \hspace{0.25cm} t \in [a,\rho(b)) \setminus \J \subset [a,\rho(b)) \\ & \\
|b - b|,              & \emph{if} \hspace{0.25cm} t \in [\rho(b),b].
\end{array}\right\} \emph{$<$} \hspace{0.25cm} \delta
\end{equation*}
 

\noindent where $\tau \in \T$ it's agreed to be a discrete point to the right of $\T$ such that $a \leq \tau < \rho(b)$ and $t \in [t_{i-1 },t_i)$ for some $t_{i-1},t_i \in P_\delta$, by the previous statement. 
 

To show this inequality, we consider the three possibilities: $t \in [a,\rho(b))$, $t \in [a,\rho(b)) \setminus \J$ and $t \in [\rho(b),b]$.

\begin{enumerate}

\item[] (i) Let $t \in [\tau,\sigma(\tau)) \subset \J$.






Then $t \in [a,\rho(b))$, where there are $t_{i-1},t_i \in P_\delta$ such that $t \in [t_{i-1},t_i)$ . Let us verify that $t_{i-1} \leq \tau \leq t <\sigma(\tau) \leq t_i$. In effect, it is enough to show that $t_{i-1} \leq \tau$ and $\sigma(\tau) \leq t_i$. 

Thus, suppose by contradiction that $\tau < t_{i-1}$. Therefore, $\sigma(\tau) \leq t_{i-1}$, then $\tau \leq t < \sigma(\tau) \leq t_{i-1} \leq t < t_i$, whence $t < t$, contradiction. Therefore $t_{i-1} \leq \tau$.

Now suppose by contradiction that $t_i < \sigma(\tau)$. Therefore $t_i \leq \tau$, then $t_{i-1} \leq t < t_i \leq\tau \leq t < \sigma(\tau)$, where $t < t$, contradiction. Therefore $\sigma(\tau) \leq t_i$.

Combining the previous results, we show that $t_{i-1} \leq \tau \leq t < \sigma(\tau) \leq t_i$.

Let us now consider the two possibilities for the pair $t_{i-1},t_i \in P_\delta$:

\begin{enumerate}

\item[•] \hspace{0.5cm} If $|t_i - t_{i-1}| < \delta$:

Then $|\sigmadelta(t) - \bigsigma(t)| = |t_i - \sigma(\tau)| \leq |t_i - t_{i-1}|< \delta$.

\item[•] \hspace{0.5cm} If $|t_i - t_{i-1}| \geq \delta$:


Then, since $P_\delta$ is delta-partition, it follows that $\sigma(t_{i-1}) = t_i$. Thus, $t_{i-1} \leq \tau \leq t \leq \sigma(\tau) \leq \sigma(t_{i-1})$, where $t_{i-1} = \tau $. Therefore, $|\sigmadelta(t) - \bigsigma(t)| = |t_i - \sigma(t_{i-1})| = 0 < \delta$.

\end{enumerate}


Therefore, if $t \in \J$, then $|\sigmadelta(t) - \bigsigma(t)| < \delta$.

\item[] (ii) Let $t \in [a,\rho(b)) \setminus \J$.


Similar to the previous item, there are $t_{i-1},t_i \in P_\delta$ such that $t \in [t_{i-1},t_i)$. Let us consider the two possibilities for the pair $t_{i-1},t_i \in P_\delta$:

\begin{enumerate}

\item[•] \hspace{0.5cm} let $|t_i - t_{i-1}| < \delta$:

Then $|\sigmadelta(t) - \bigsigma(t)| = |t_i - t| \leq |t_i - t_{i-1}|< \delta$.

\item[•] \hspace{0.5cm} Let $|t_i - t_{i-1}| \geq \delta$:


Since $P_\delta$ is delta-partition, then it follows that $\sigma(t_{i-1}) = t_i$. Thus, $t \in [t_{i-1},\sigma(t_{i-1}))$ where $t_{i-1} \in [a,\rho(b))$ is a discrete point to the right of $\T$. Therefore $t \in \J$, absurd, because $t \in [a,\rho(b)) \setminus \J$. Thus, if $t \in [a,\rho(b)) \setminus \J$, then $|t_i - t_{i-1}| < \delta$, where $|\sigmadelta(t) - \bigsigma(t)| < \delta$.

\end{enumerate}

Therefore, if $t \in [a,\rho(b)) \setminus \J$, then $|\sigmadelta(t) - \bigsigma(t)| < \delta$.

\item[] (iii) Let $t \in [\rho(b),b]$.


If $b$ is left dense, then $t = b$. Therefore $|\sigmadelta(b) - \bigsigma(b)| = |b - b| = 0 < \delta$. If $b$ is right discrete, then $t_{n-1} \in P_\delta$ is such that $t_{n-1} = \rho(b)$, since $\delta < \delta_0 $. Therefore $\sigmadelta(t) = b$. Therefore, $|\sigmadelta(t) - \bigsigma(t)| = |b - b| = 0 < \delta$.

\end{enumerate} 


Thus, for all $\delta \in (0,\delta_0)$ we have $|\sigmadelta(t) - \bigsigma(t)| < \delta$ for all $t \in [a,b]$ and independent of the partition $P_\delta$. This shows that $\sigmadelta$ converges uniformly to $\bigsigma$ in $[a,b]$ in the sense that we defined in the statement of this theorem. \hfill $\square$

\end{prova}


\begin{proposicao}

The function $\bigsigma:[a,b] \longrightarrow \R$ defined in \emph{Proposition \ref{prop_bigsigma}} has the following properties
\begin{enumerate}

\item[]\textbf{\emph{(i)}} $S(\T)(t) = \sigma_\T(t)$ for all $t \in [a,\rho(b))_\T$.

\item[]\textbf{\emph{(ii)}} $t \leq S(\T)(t) \leq t \hspace{0.15cm} + \hspace{-0.1cm} \disp\sup_{\tau \in [a,\rho(b))_\T} \mu_\T(t)$ for all $t \in [a,b]$.

\item[]\textbf{\emph{(iii)}} $S(\T)$ is non-decreasing.

\end{enumerate}

\label{prop_?}
\end{proposicao}

\begin{prova}



We will denote, for simplicity of notations, $\sigma_\T = \sigma$, $\rho_\T = \rho$ and $\mu_\T = \mu$. 

(i) Once the function $\bigsigma:[a,b] \longrightarrow \R$ has been defined according to the previous proposition, we can see that $\bigsigma$ is actually an extension of the jump operator $\sigma: \T \longrightarrow \T$ to $[a,\rho(b))$, that is, $\bigsigma(t) = \sigma(t)$ for all $t \in [a,\rho(b))_\T $.

Indeed, if $t \in [a,\rho(b))_\T$ is right-discrete in $\T$, then $t \in [t,\sigma(t)) \subset \J$ , where $\bigsigma(t) = \sigma(t)$. On the other hand, if $t \in [a,\rho(b))_\T$ is right-dense in $\T$, then $t \notin \J$, since if $t \in \J $, then there would be right-discrete $\tau \in \T$ such that $t \in [\tau, \sigma(\tau))$, where $t = \tau$ and $t$ would be right-discrete. Therefore $t \in [a,\rho(b))\setminus \J$ and therefore $\bigsigma(t) = t = \sigma(t)$.



(ii) Let $t \in [a,\rho(b))$. If $t \in [\tau, \sigma(\tau)) \subset \J$, then $t \leq \sigma(\tau) = \bigsigma(t)$. Furthermore, if $t \in [a,\rho(b)) \setminus \J$, then $t \leq t = \bigsigma(t)$. Therefore $t \leq \bigsigma(t)$ for all $t \in [a,\rho(b))$. If $t \in [\rho(b),b]$, then $t \leq b = \bigsigma(t)$. This shows the first inequality.

For the second inequality, let $t \in [a,\rho(b))$ again. If $t \in [\tau, \sigma(\tau)) \subset \J$, then $\bigsigma(t) = \sigma(\tau) = \tau + \mu(\tau) \leq \tau + \sup \{\mu(\tau) :\tau \in [a,\rho(b))_\T \} \leq t + \sup \{\mu(\tau) :\tau \in [ a,\rho(b))_\T \}$, since $\tau \in [a,\rho(b))_\T$ and $\tau \leq t$. If $t \in [a,\rho(b)) \setminus \J$, then $\bigsigma(t) = t \leq t +\sup \{\mu(\tau) :\tau \in [a ,\rho(b))_\T \}$, since $\mu: \T \to \T$ is a non-negative function. Furthermore, if $t \in [\rho(b),b]$, then $\bigsigma(t) = b \leq b + \sup \{\mu(\tau) :\tau \in [a,\rho(b))_\T \}$.






(iii) Let $x,y \in [a,\rho(b))$, with $x \leq y$. Let's show that $\bigsigma(x) \leq \bigsigma(y)$. Consider the cases:

(iii).(i) If $x,y \in \J$, then there exist right-discrete $\tau_1,\tau_2 \in [a,\rho(b))_\T$ such that $x \in [\tau_1,\sigma(\tau_1))$ and $y \in [\tau_2,\sigma(\tau_2))$. Suppose by contradiction that $\tau_1 > \tau_2$. Thus, we have $\tau_2 < \tau_1 \leq x \leq y < \sigma(\tau_2)$, where $\tau_2 < \tau_1 < \sigma(\tau_2)$, contradiction. Therefore $ \tau_1 \leq \tau_2$, where $\bigsigma(x) = \sigma(\tau_1) \leq \sigma(\tau_2) = \bigsigma(y)$;

(iii).(ii) If $x,y \in [a,\rho(b)) \setminus \J$, then $\bigsigma(x) = x \leq y = \bigsigma(y)$;

(iii).(iii) If $x \in \J$ and $y \in [a,\rho(b)) \setminus \J$, then there is $\tau_1 \in [a,\rho(b) )_\T$ right-discrete such that $x \in [\tau_1,\sigma(\tau_1))$. See that $y \notin (-\infty, \tau_1)$, since $x \leq y$, and $y \notin [\tau_1, \sigma(\tau_1))$, since otherwise we would have $y \in \J$. Therefore $y \geq \sigma(\tau_1)$, where $\bigsigma(y) = y \geq \sigma(\tau_1) = \bigsigma(x)$;

(iii).(iv) If $x \in [a,\rho(b)) \setminus \J$ and $y \in \J$, then there exists right-discrete $\tau_2 \in \T$ such that $y \in [\tau_2, \sigma(\tau_2))$. Therefore $\bigsigma(x) = x \leq y < \sigma(\tau_2) = \bigsigma(y)$.


Therefore for all $x,y \in [a,\rho(b))$, with $x \leq y$, we have $\bigsigma(x) \leq \bigsigma(y)$. Now, if $x,y \in [\rho(b), b]$, then $\bigsigma(x) = b \leq b = \bigsigma (y)$. Finally, if $x \in [a,\rho(b))$ and $y \in [\rho(b),b]$, then $\bigsigma(x) = \sigma(\tau_1)$ for some right-discrete $\tau_1 \in [a,\rho(b))_\T$. Note that $\tau_1 < \rho(b)$ implies that $\sigma(\tau_1) \leq \rho(b)$. Therefore $\bigsigma(x) = \sigma(\tau_1) \leq \rho(b) \leq b = \bigsigma(y)$. This concludes result (iii). \hfill $\square$
\end{prova}



Next, we state Abel's Summation Lemma, which is the main tool for demonstrating the Theorem \ref{teo_identidade}.

\newpage

\begin{lema}

\emph{(Abel's Summation)}. Let $(\alpha_i)_{i \in \N}, (\beta_i)_{i \in \N}$ be sequences of complex numbers. For all $n \in \N$, with $n \geq 2$, equality is verified
\begin{equation*}
\disp\sum_{i = 1}^{n} \alpha_i \cdot \beta_i = \disp\sum_{i=1}^{n-1} \left( \disp\sum_{j=1}^{i} \alpha_j \right) \cdot (\beta_i - \beta_{i+1}) + \beta_n \cdot \disp\sum_{i=1}^{n} \alpha_i 
\end{equation*}

\label{lema_abel}
\end{lema}

\begin{prova}
See Theorem 4.2. from \cite{introduction}. 
\end{prova}


\section{Main Results}

\begin{teorema}

Let $a,b \in \T$, with $a < \rho_\T(b)$, and $\T,\Ttil$ be time scales. Then there exists $S:[a,b]_{\Ttil} \longrightarrow \R$ with the following properties:
\begin{enumerate}

\item[]\textbf{\emph{(i)}} $S \defi \bigsigma_{|_{\Ttil}}$ is non-decreasing, where $\bigsigma$ is given as in \emph{Proposition \ref{prop_bigsigma}}.

\item[]\textbf{\emph{(ii)}} $S(t) = \sigma_\T(t)$ for all $t \in [a,\rho(b))_\T$.

\item[]\textbf{\emph{(iii)}} $a \leq t \leq S(t) \leq t \hspace{0.1cm} + \hspace{-0.25cm} 
\disp\sup_{s \hspace{0.05cm} \in \hspace{0.05cm} [a,\rho_\T(b))_\T} \hspace{-0.25cm} \mu_\T(s) \leq b \hspace{0.25cm}$ for all $t \in [a,\rho_\T(b))_{\Ttil}$.

\item[]\textbf{\emph{(iv)}} For all $f \in $ \emph{C}$_{\emph{rd}}^1([a,b]_{\Ttil},\R)$ we have the following $\Delta$-Riemann integrals:

\begin{equation}
\displaystyle\int_{a}^{b}  f(s) \Delta_\T s \; = \; [s \cdot f(s)]\Big|_{a}^b - \displaystyle\int_{a}^{b} (S\cdot f^{\tilde{\vartriangle}})(s) \Delta_{\Ttil} s.
\end{equation}

\end{enumerate}


\noindent where $\rho_\T,\sigma_\T: \T \longrightarrow \T$ are respectively the backward and forward operators of $\T$, $\mu_\T: \T \to \R^+ _0$ is the graininess operator of $\T$ and $f^{\tilde{\deri}}$ is the delta-derivative of $f$ on the scale $\Ttil$. Nonetheless, $S$ can be taken as $$S(t) = \inf\{s \in \T: s > t\}$$ 
for all $t \in [a,b]_{\Ttil}$.

\label{teo_identidade}
\end{teorema}

\begin{prova}


Let $a,b \in \T \subset \Ttil$ with $a < \rho_\T(b)$. Take $S:[a,b]_{\Ttil} \longrightarrow \R$ for $S = \bigsigma_{|_{\Ttil}}$. By Proposition \ref{prop_bigsigma}, we immediately verify items (i), (ii) and (iii). Now we will verify that item (iv) is also true. Thus, take the application $f \in $ \emph{\emph{C}}$_{\emph{\emph{rd}}}^1([a,b]_{\Ttil},\R)$.

First, let's see that the identity in (iv) makes sense. Since $f: [a,b]_{\Ttil} \longrightarrow \mathbb{R}$ is continuous, then $f_{|_\mathbb{T}}:[a,b]_{\T} \longrightarrow \mathbb{R}$, your restriction to $[a,b]_{\T}$ (a closed subset of your domain) is continuous on $[a,b]_{\T}$, and therefore $ \Delta$-Riemann integrable.


Therefore, let $P_\delta \in \mathcal{P}_\delta(a,b)$, with $P_\delta = \{a = t_1, t_2, \cdots, t_{n-1}, t_n = b\}$, a delta-partition of the interval $[a,b]_\T$. We take the $S_\delta$ $\Delta$-sum of \emph{Riemann} for $f$ for the partition $P_\delta$. Then $S_\delta$ is given by
\begin{equation*} 
S_\delta = \disp\sum_{i=2}^{n} f(\xi_i)(t_i - t_{i-1}).
\end{equation*}


\noindent where $\xi_i \in [t_{i-1},t_i)$ is chosen arbitrarily. As $f$ is integrable, it follows that the limit of $S_\delta$ when $\delta$ tends to $0^+$ exists and is independent of the choices of $P_\delta$ and $\xi_i$. This limit is defined as the integral of $f$ in $[a,b]_\mathbb{T}$. Thus, for each $P_\delta \in \mathcal{P}(a,b)$ let us fix $\xi_i = t_{i-1} \in P_\delta$, with $i \in \{2, \cdots, n\}$, where $S_\delta$ is given by
\begin{equation}\label{S_definition}
S_\delta = \disp\sum_{i=2}^{n} f(t_{i-1})(t_i - t_{i-1}).
\end{equation}


Using Abel's summation lemma, we can rewrite the identity (\ref{S_definition}) in the form
\begin{equation} \label{S_second_definition}
S_\delta = \disp\sum_{i=2}^{n-1} \left( \disp\sum_{j=2}^{i}
(t_j - t_{j-1}) \right) \cdot (f(t_{i-1})-f(t_i)) + \left(\disp\sum_{i=2}^{n} (t_i - t_{i-1}) \right) \cdot f(t_{n-1}). 
\end{equation}

\newpage


Since $f$ is differentiable in 
$[a,b]_{\Ttil}$ and $P_\delta \subset [a,b]_{\Ttil}$, taking $f^{\tilde{\vartriangle}}$ as the derivative of $f$ in $[a,b]_{\Ttil}$, we have
\begin{equation}
f(t_{i-1}) - f(t_i) = - \disp\int_{t_{i-1}}^{t_i} f^{\tilde{\vartriangle}}(s) \Delta_{\Ttil} s.
\end{equation}


Furthermore, for all $i \in \{2, \cdots, n-1\}$, we have
\begin{eqnarray*}
\disp\sum_{j=2}^{i} (t_j - t_{j-1}) & = &(t_2 - t_1) + (t_3 - t_2) + \cdots + (t_i - t_{i-1}) =  t_i - a. 
\end{eqnarray*}


We can rewrite equation (\ref{S_second_definition}) in the new form
\begin{equation}
S_\delta = - \disp\sum_{i=2}^{n-1} (t_i - a) \cdot \disp\int_{t_{i-1}}^{t_i} f^{\tilde{\vartriangle}}(s) \Delta_{\Ttil} s + (b-a)f(t_{n-1}). 
\end{equation}

Then 
\begin{eqnarray*}
S_\delta & = & - \disp\sum_{i=2}^{n-1} t_i \cdot \disp\int_{t_{i-1}}^{t_i} f^{\tilde{\vartriangle}}(s) \Delta_{\Ttil} s \hspace{0.1cm} + \hspace{0.1cm} a \disp\sum_{i=2}^{n-1} \disp\int_{t_{i-1}}^{t_i} f^{\tilde{\vartriangle}}(s) \tilde{\Delta} s \hspace{0.1cm} + \hspace{0.1cm} (b-a)f(t_{n-1}) \\
& = & - \disp\sum_{i=2}^{n-1} \disp\int_{t_{i-1}}^{t_i}  t_i \cdot f^{\tilde{\vartriangle}}(s) \Delta_{\Ttil} s \hspace{0.1cm} + \hspace{0.1cm}  a \disp\int_{a}^{t_{n-1}} f^{\tilde{\vartriangle}}(s) \tilde{\Delta} s \hspace{0.1cm} + \hspace{0.1cm} (b-a)f(t_{n-1}) \\
& = & - \disp\sum_{i=2}^{n-1}\disp\int_{t_{i-1}}^{t_i} t_i \cdot f^{\tilde{\vartriangle}}(s) \Delta_{\Ttil} s \hspace{0.1cm} + \hspace{0.1cm} a(f(t_{n-1})- f(a)) \hspace{0.1cm} + \hspace{0.1cm} (b-a)f(t_{n-1}). 
\end{eqnarray*}


As $S_{P_\delta}$ is integrable and it holds that $S_{P_\delta}(s) = t_i$ for $ s \in [t_{i-1},t_i)$, with $i \in \{2, \cdots, n-1\}$ and $S_{P_\delta}(s) = b$ for $ s \in [t_{n-1},b]$, we rewrite the above identity in the form
\begin{eqnarray*}
S_\delta & = & - \disp\sum_{i=2}^{n-1}\disp\int_{t_{i-1}}^{t_i} S_{P_\delta}(s) \cdot f^{\tilde{\vartriangle}}(s) \Delta_{\Ttil} s \hspace{0.1cm} + \hspace{0.1cm} a(f(t_{n-1})- f(a)) \hspace{0.1cm} + \hspace{0.1cm} (b-a)f(t_{n-1}). \\
& = & - \disp\int_{a}^{t_{n-1}} (S_{P_\delta} \cdot f^{\tilde{\vartriangle}})(s) \Delta_{\Ttil}s \hspace{0.1cm} + \hspace{0.1cm} a(f(t_{n-1})- f(a)) \hspace{0.1cm} + \hspace{0.1cm} (b-a)f(t_{n-1})\\
& = & - \disp\int_{a}^{t_{n-1}} (S_{P_\delta}\cdot f^{\tilde{\vartriangle}})(s) \Delta_{\Ttil} s \hspace{0.1cm} - \hspace{0.1cm}  \disp\int_{t_{n-1}}^{b} (S_{P_\delta} \cdot f^{\tilde{\vartriangle}})(s) \Delta_{\Ttil} s \hspace{0.1cm} + \\ & &  + \disp\int_{t_{n-1}}^{b} (S_{P_\delta} \cdot f^{\tilde{\vartriangle}})(s) \Delta_{\Ttil} s \hspace{0.1cm} + \hspace{0.1cm}  a(f(t_{n-1})- f(a)) \hspace{0.1cm} + \hspace{0.1cm} (b-a)f(t_{n-1}) \\
& = & - \disp\int_{a}^{b} (S_{P_\delta} \cdot f^{\tilde{\vartriangle}})(s) \Delta_{\Ttil} s \hspace{0.1cm} + \hspace{0.1cm} b(f(b) - f(t_{n-1})) \hspace{0.1cm} + \hspace{0.1cm}  a(f(t_{n-1})- f(a)) \hspace{0.1cm} + \hspace{0.1cm} (b-a)f(t_{n-1}) \\
& = & - \disp\int_{a}^{b} (S_{P_\delta} \cdot f^{\tilde{\vartriangle}})(s) \Delta_{\Ttil} s  
\hspace{0.1cm} + \hspace{0.1cm} bf(b) - af(a). 
\end{eqnarray*}


Thus, taking $\disp\lim_{\delta \to 0^+} S_\delta$, it follows that
\begin{eqnarray*}
\disp\lim_{\delta \to 0^+} S_\delta & = & (bf(b) -af(a))\hspace{0.1cm} - \hspace{0.1cm} \disp\lim_{\delta \to 0^+} \disp\int_{a}^{b} (S_{P_\delta}  \cdot f^{\tilde{\vartriangle}})(s) \Delta_{\Ttil} s \label{a}\\ 
& \overset{(a)}{=}& (bf(b) -af(a)) \hspace{0.1cm} - \hspace{0.1cm} \disp\int_{a}^{b}  \disp\lim_{\delta \to 0^+} S_{P_\delta}(s) \cdot f^{\tilde{\vartriangle}}(s) \Delta_{\Ttil} s \label{b}\\ 
& \overset{(b)}{=} & (bf(b) -af(a)) \hspace{0.1cm} - \hspace{0.1cm} \disp\int_{a}^{b}S(s) \cdot f^{\tilde{\vartriangle}}(s) \Delta_{\Ttil} s. 
\end{eqnarray*}



where passages (a) and (b) are true because, given $\varepsilon > 0$ and $(\delta_n)_{n \in \N}$ a positive sequence that converges to $0$, by uniform convergence given by the Theorem \ref{prop_?} from $\sigmadelta$ to $\bigsigma$ in $[a,b]$, we have that there is $\delta_0 > 0$ such that $|\sigmadelta(t) - \bigsigma(t )| < \varepsilon$ for all $\delta \in (0,\delta_0)$, for all $P_\delta \in \mathcal{P}(a,b)$ and for all $t \in [a,b]$.

Therefore $|\mbox{S}(t) - \mbox{S}_{P_\delta}(t)| < \varepsilon$ for all $\delta \in (0,\delta_0)$, for all $P_\delta \in \mathcal{P}(a,b)$ and for all $t \in [a,b] _{\Ttil}$. Thus, take $n_0 \in \N$ such that $1/n_0 < \delta_0$. Therefore, $(\delta_n)_{n \geq n_0} \subset (0, \delta_0)$. Furthermore, as $f^{\tilde{\deri}} \in \mbox{C}_{\emph{\emph{rd}}}([a,b]_{\Ttil},\R)$, then $f^{\tilde{\deri}}$ is bounded in $[a,b]_{\Ttil}$. Take $M > 0$ such that $|f^{\tilde{\deri}}(t)| \leq M$ for all $t \in [a,b]_{\Ttil}$. That way,
\begin{eqnarray*}
\left| \disp\int_{a}^{b} \mbox{S}_{P_{\delta_n}}(s)f^{\tilde{\deri}}(s) \Delta_{\Ttil} s  - \disp\int_{a}^{b}\mbox{S}(s)f^{\tilde{\deri}}(s) \Delta_{\Ttil} s  \right| & = & \left| \disp\int_{a}^{b} (\mbox{S}_{P_{\delta_n}}(s) - \mbox{S}(s))  f^{\tilde{\deri}}(s) \Delta_{\Ttil} s  \right|  \\
& \leq & \disp\int_{a}^{b} \left|  \mbox{S}_{P_{\delta_n}}(s) - \mbox{S}(s) \right| \cdot \left| f^{\tilde{\deri}}(s) \right|  \Delta_{\Ttil} s  \\
& < & M  \disp\int_{a}^{b} \varepsilon \;\; \Delta_{\Ttil} s \; = \; M(b-a) \varepsilon.
\end{eqnarray*}

This shows that $\lim_{n \in \N} \smallint_{a}^{b} \mbox{S}_{P_{\delta_n}}(s)f^{\tilde{\deri}}(s) \tilde{\Delta} s = \smallint_{a}^{b}\mbox{S}(s)f^{\tilde{\deri}}(s) \tilde{\Delta} s$. As $(\delta_n)_{n \in \N}$ is an arbitrary sequence, thus $\lim_{\delta \to 0+} \smallint_{a}^{b} \mbox{S}_{P_{\delta}}(s)f^{\tilde{\deri}}(s) \tilde{\Delta} s = \smallint_{a}^{b}\mbox{S}(s)f^{\tilde{\deri}}(s) \tilde{\Delta} s$. Therefore,
\begin{equation}
\disp\int\limits_{a}^{b} f(s) \Delta_\T s  =  \disp\lim_{\delta \to 0^+} S_\delta =  bf(b) -af(a) \hspace{0.1cm} - \hspace{0.1cm} \disp\int_{a}^{b} S(s) \cdot f^{\tilde{\vartriangle}}(s) \Delta_{\Ttil} s
\end{equation}

\hfill $\square$
\end{prova}


\begin{exemplo} 

(Integration by Parts) Let $a,b \in \T$ and $f \in$ be $C_{\emph{rd}}(\T,\R)$. Let us take $\Ttil = \T$. Note that $\bigsigma_\T(t) = \sigma_\T(t)$ for all $t \in [a,b)_\T$. Therefore, by the \emph{Theorem \ref{teo_identidade}}, we have
\begin{equation*}
\displaystyle\int_{a}^{b}  f(s) \Delta_\T s = [s \cdot f(s)]\Big|_{a}^b - \displaystyle\int_{a}^{b} (\sigma_\T \cdot f^{\vartriangle})(s) \Delta_\T s.
\end{equation*}


We just need to note that:
\begin{equation*}
\displaystyle\int_{\rho_\T(b)}^{b} \sigma_\T(s) \cdot f^{\vartriangle}(s) \Delta_\T s = \displaystyle\int_{\rho_\T(b)}^{b} (\bigsigma_{\mathbb{T}} \cdot f^{\vartriangle})(s) \Delta_\T s.
\end{equation*}


Indeed, if $b$ is left dense then the above identity is trivially satisfied. If $b$ is left-discrete, then $\rho_\T(b)$ is right-discrete, so $\sigma_\T(\rho_\T(b)) = b = \bigsigma_\T(\rho_\T (b))$. This way we have to
\begin{equation*}
\displaystyle\int_{\rho_\T(b)}^{b} \sigma(s) \cdot f^{\vartriangle}(s) \Delta_\T s = \mu_\T(\rho_\T(b)) \cdot \sigma_\T(\rho_\T(b)) \cdot f^{\vartriangle}(\rho_\T(b)) = \displaystyle\int_{\rho_\T(b)}^{b} (\bigsigma_{\mathbb{T}} \cdot f^{\vartriangle})(s) \Delta_\T s.
\end{equation*}

Thus
\begin{equation*}
\displaystyle\int_{a}^{b} \sigma(s) \cdot f^{\vartriangle}(s) \Delta s = \displaystyle\int_{a}^{b} (\bigsigma_{\mathbb{T}} \cdot f^{\vartriangle})(s) \Delta s.
\end{equation*}

\end{exemplo}


\begin{exemplo} 


Let $a,b \in \N \subset \R$ and $f \in$ be $C_{\emph{rd}}^1(\R,\R)$. Let us take $\T = \N$ and $\Ttil = \R$. Let's see that, for $\N$, we have $\bigsigma(t) = \lfloor t \rfloor + 1$ for all $t \in [a,b)$ and $\bigsigma(b) = b$. In effect, let $t \in [a,b)$. As we saw previously, the set $\J$ defined for the interval $[a,b]_\N$ is given by $I(\N) = [a,b-1)$. See that $[a,\rho(b)) \setminus I(\N) = \emptyset$.

In general we have that $t \in [\lfloor t\rfloor,\lfloor t \rfloor+1)$, where $\lfloor t\rfloor,\lfloor t \rfloor+1 \in \N$. Thus, if $t \in [a,b-1)$ then $S(\N)(t) = \lfloor t \rfloor + 1$, if $t \in [b-1,b]$, then $S(\N)(t) = b = \lfloor t \rfloor + 1$. Note that $\lfloor t + 1 \rfloor= \lfloor t \rfloor + 1$ for all $t \in \R$. 

It follows from the Theorem \ref{teo_identidade} that
\begin{eqnarray*}
\displaystyle\int_{a}^{b}  f(s) \Delta_{\N} s & = & [s \cdot f(s)]\Big|_{a}^b - \displaystyle\int_{a}^{b} (S(\N) \cdot f')(s) \, \Delta_{\R} s. \\
\end{eqnarray*}

\newpage


Simplifying these expressions, we have
\begin{eqnarray*}
\displaystyle\int_{a}^{b}  f(s) \Delta_{\N} s & = & \displaystyle\sum_{i = a}^{b-1}  f(i), \quad \disp\int_{a}^{b} (S(\N) \cdot f')(s) \, \Delta_{\R} s \; = \; \disp\int_{a}^{b} (S(\N) \cdot f')(s) \, d s  \\
\Rightarrow \displaystyle\sum_{i = a}^{b-1}  f(i) & = & [s \cdot f(s)]\Big|_{a}^b - \displaystyle\int_{a}^{b} (\lfloor s \rfloor + 1) \cdot f'(s)  \hspace{0.1cm} ds \\
& = & (b-1)f(b) - (a-1)f(a) - \displaystyle\int_{a}^{b} \lfloor s \rfloor \cdot f'(s)  \hspace{0.1cm} ds. 
\end{eqnarray*}

\end{exemplo}


\vspace{0.15cm}

\begin{corolario}

Let $a,b \in \T \subset \Ttil, a \leq b$ and $f \in$ \emph{C}$_{\emph{rd}}^1(\Ttil,\R)$. We have that, if $f$ is non-decreasing in $[a,b]_{\Ttil}$, then
\begin{equation*}
\displaystyle\int_{a}^{b} f(s) \Delta_\T s \leq \displaystyle\int_{a}^{b} f(s) \Delta_{\Ttil} s.
\end{equation*}

\end{corolario}

\begin{prova}

Note that, if $f$ is non-decreasing in $[a,b]_{\Ttil}$, then $f^{\tilde{\vartriangle}} \geq 0$ in $[a,b)_{\Ttil}$. Analogously, if $f$ is non-increasing in $[a,b]_{\Ttil}$ then $f^{\tilde{\vartriangle}} \leq 0$ in $[a,b)_{\Ttil}$. First, let us verify that 
$\bigsigma(s) \geq \sigma_{\Ttil}(s)$ in $[a,b)_{\Ttil}$. Let $t \in [a,b)_{\Ttil}$. If $t \in [\tau, \sigma_\T(\tau)) \subset \J$, then $\bigsigma(t) = \sigma_\T(\tau)$. As 
$t,\sigma_\T(\tau) \in \Ttil$ and $ t < \sigma_\T(\tau) \Rightarrow \sigma_{\Ttil}(t) \leq \sigma_\T(\tau)$, whence $\bigsigma(t) \geq \sigma_{\Ttil}(t)$. If $t \notin \J$, then 
$\bigsigma(t) = t \leq \sigma_{\Ttil}(t)$.

By the theorem \ref{teo_identidade} we have the identity
\begin{equation*}
\displaystyle\int_{a}^{b}  f(s) \Delta_\T s = [s \cdot f(s)]\Big|_{a}^b - \displaystyle\int_{a}^{b} (\bigsigma_{\tilde{\mathbb{T}}} \cdot f^{\tilde{\vartriangle}})(s) \Delta_{\Ttil} s.
\end{equation*}


Therefore, if $f$ is non-decreasing in $[a,b)_{\Ttil}$ then $\bigsigma_{\Ttil}(s) \cdot f^{\tilde{\vartriangle}}(s) \geq \sigma_{\Ttil}(s) \cdot f^{\tilde{\vartriangle}}(s)$, where we have that
\begin{eqnarray*}
\displaystyle\int_{a}^{b}  f_{\mathbb{T}}(s) \Delta_\T s & = & [s \cdot f(s)]\Big|_{a}^b - \displaystyle\int_{a}^{b} (\bigsigma_{\tilde{\mathbb{T}}} \cdot f^{\tilde{\vartriangle}})(s) \Delta_{\Ttil} s \\
 & \leq & [s \cdot f(s)]\Big|_{a}^b - \displaystyle\int_{a}^{b} (\sigma_{\Ttil} \cdot f^{\tilde{\vartriangle}})(s) \Delta_{\Ttil} s \; = \; \displaystyle\int_{a}^{b} f(s) \Delta_{\Ttil} s.
\end{eqnarray*}

\hfill $\square$
\end{prova}


\vspace{0.15cm}

\begin{teorema}

Let $\{\T_n\}_{n \in \N_0}$ be a family of time scales such that $a,b \in \T_0 \subset \T_1 \subset ... \subset \bigcup_{n \in \ N_0} \T_n = \Ttil$, with
$\inf \T_0 < a < \rho_{\T_0}(b) \leq b < \sup \T_0$. If $f \in$ \emph{C}$_{\emph{rd}}^1([a,b]_{\Ttil},\R)$ and 
$\bigsigman$ converges punctually to $\bigsigmatil$ in $[a,b]$, then
\begin{equation*}
\disp\int_{a}^{b} f(s) \Delta_{\Ttil} s \; = \; \disp\lim_{n \to \infty} \disp\int_{a}^{b} f(s) \Delta_{\T_n} s.
\end{equation*}

\end{teorema}

\begin{prova}

Since $f$ is continuous on $[a,b]_{\Ttil}$, as $\T_n \subset \Ttil$ are all closed sets, then the restriction of $f$ to $[a,b] _{\T_n}$, denoted by $f_{|_{\T_n}}$, continue to be continuous and therefore integrable. Thus, let $M > 0$ be such that $f^{\tilde{\deri}}(s) \leq M$ for all $s \in [a,b]_{\Ttil}$. For all $n \in \N$ we have that $a,b \in \T_n \subset \Ttil$, where $a < \rho_{\T_0}(b) \leq \rho_{\T_n}(b) $, since $\T_0 \subset \T_n$ and $\rho_{\T_0}(b)> a > \inf \T_0$. We can then apply the Theorem \ref{teo_identidade} for each pair $\T_n, \Ttil$, where
\begin{equation*}
\disp\int_{a}^{b} f(s) \Delta_{\T_n} s  \; = \;  [s\cdot f(s)]\Big|_a^{b} - \displaystyle\int_{a}^{b} (S(\T_n)_{\Ttil} \cdot f^{\tilde{\vartriangle}})(s) \Delta_{\Ttil} s.
\end{equation*}


As $f$ is $\Delta_{\Ttil}$-differentiable, by the integration by parts method, we have
\begin{equation*}
\disp\int_{a}^{b} f(s) \Delta_{\Ttil} s  \; = \;  [s\cdot f(s)]\Big|_a^{b} - \displaystyle\int_{a}^{b} (\sigma_{\Ttil} \cdot f^{\tilde{\vartriangle}})(s) \Delta_{\Ttil} s.
\end{equation*}

\newpage


Furthermore, as $\bigsigmatil_{|_{\Ttil}}$ is an extension of $\sigma_{\Ttil}$ in $[a,\rho_{\Ttil}(b))$. This way, we can rewrite the above identity in the form
\begin{equation*}
\disp\int_{a}^{b} f(s) \Delta_{\Ttil} s  \; = \;  [s\cdot f(s)]\Big|_a^{b} - \displaystyle\int_{a}^{b} (S(\Ttil) \cdot f^{\tilde{\vartriangle}})(s) \Delta_{\Ttil} s.
\end{equation*}


Comparing equations, we have for each $n \in \N_0$
\begin{eqnarray*}
\Big| \disp\int_{a}^{b} f(s) \Delta_{\T_n} s - \disp\int_{a}^{b} f(s) \Delta_{\Ttil} s \Big| 
& = & \Big| \disp\int_{a}^{b} f^{\tilde{\vartriangle}} \cdot (S(\Ttil) - S(\T_n)_{\Ttil}) (s) \Delta_{\Ttil} s \Big| \\
& \leq & M \disp\int_{a}^{b} |S(\Ttil) - S(\T_n)_{\Ttil} |(s) \Delta_{\Ttil} s. 
\end{eqnarray*}


Since $\{\T_n\}_{n \in \N}$ is an ascending chain of scales and $\inf \T_{n+1} \leq \inf \T_n < a < \rho_0(b) \leq \rho_n (b) \leq \rho_{n+1}(b) \leq b < \sup \T_n \leq \sup \T_{n+1}$ for all $n \in \N_0$, it follows that $ \bigsigmatil(t) \leq ... \leq S(\T_1)(t) \leq S(\T_0)(t)$ for all $t \in [a,b]$, since $\T \subset \Ttil$ implies that $S(\Ttil)(t) \leq S(\T)(t)$ for all $t \in [a,b]$. Note that as $S(\T_n):[a,b]_{\Ttil} \longrightarrow \R$ converges pointwise to $S(\Ttil):[a,b]_{\Ttil} \longrightarrow \R $, with $[a,b]_{\Ttil}$ compact and $S(\T_n)$ a monotone sequence. Then, by Dini's Theorem (\cite{rudin}), it follows that $S(\T_n)$ converges uniformly to $S(\Ttil)$, which concludes the result.

\hfill $\square$
\end{prova}


\begin{corolario}

Let $a,b \in \T$, with $a \leq b$, and $\T \subset \Ttil$ be time scales. Then, for every $\Delta_{\Ttil}$-integrable function $f:[a,b]_{\Ttil} \longrightarrow \R$  the following equality between the $\Delta$-Riemann integrals holds
\begin{eqnarray}
\int_{a}^{b} f(s) \, \Delta_\T s & = & \int_{a}^{b} f \circ \sigma_{\Ttil}(s) \, \Delta_{\Ttil} s + \sum_{\tau < \sigma_\T(\tau)}^{\tau < b} \mu_\T(\tau) \cdot f(\tau) - \int_{\tau}^{\sigma_\T(\tau)} f \circ \sigma_{\Ttil}(s) \, \Delta_{\Ttil} s.
\label{eq_fim}
\end{eqnarray}

\label{cor_fim}
\end{corolario}

\begin{prova}


Suppose $a < \rho_\T(b)$ (other cases are treated trivially). Note that, as every $\Delta$-integrable function on a scale is necessarily a continuous function at almost every point, that is, it is only discontinuous on a set of zero measure, this tells us that we can demonstrate the identity (\ref {eq_fim}) only for functions $C_{rd}([a,b]_{\Ttil},\R)$. 

Given $S$ as in Theorem \ref{teo_identidade}, we have that $S(t) = t + \sum_{\tau < \sigma_\T(\tau),\, \tau < b} (\sigma_\T (\tau) - t) \cdot \chi_{[\tau, \sigma_\T(\tau))}(t)$ for all $t \in [a,b]$, where $\chi_A$ represents the characteristic function of the set $A \subset \R$. Therefore, by the Theorem \ref{teo_identidade}, for any $f \in $ C$_{rd}^1([a,b]_{\Ttil},\R)$, we have the following identity:
\begin{eqnarray}
\int_{a}^{b} f(s) \, \Delta_\T s & = & [s \cdot f(s)]\Big|_{a}^{b} - \int_{a}^{b} \left(s + \sum_{\tau < \sigma_\T(\tau)}^{\tau < b} (\sigma_\T(\tau) - s) \cdot \chi_{[\tau, \sigma_\T(\tau))}(s) \right)_{|_{\Ttil}} \hspace{-0.2cm} \cdot f^{\tilde{\vartriangle}}(s) \, \Delta_{\Ttil} s \non \\
& = & \int_{a}^{b} f \circ \sigma_{\Ttil}(s) \, \Delta_{\Ttil} s \; - \sum_{\tau < \sigma_\T(\tau)}^{\tau < b} \sigma_\T(\tau) \int_{a}^{b} \chi_{[\tau, \sigma_{\T})}(s) \cdot f^{\tilde{\vartriangle}}(s) \Delta_{\Ttil} s \non \\
&&+\; \sum_{\tau < \sigma_\T(\tau)}^{\tau < b} \int_{a}^{b} \chi_{[\tau, \sigma_{\T})}(s) \cdot s \cdot f^{\tilde{\vartriangle}}(s) \Delta_{\Ttil} s \non\\
& = & \int_{a}^{b} f \circ \sigma_{\Ttil}(s) \, \Delta_{\Ttil} s \; - \sum_{\tau < \sigma_\T(\tau)}^{\tau < b} \sigma_\T(\tau) \int_{\tau}^{\sigma_\T(\tau)}  f^{\tilde{\vartriangle}}(s) \Delta_{\Ttil} s \; + \sum_{\tau < \sigma_\T(\tau)}^{\tau < b} \int_{\tau}^{\sigma_\T(\tau)} s \cdot f^{\tilde{\vartriangle}}(s) \Delta_{\Ttil} s \non \\
& = & \int_{a}^{b} f \circ \sigma_{\Ttil}(s) \, \Delta_{\Ttil} s \; - \sum_{\tau < \sigma_\T(\tau)}^{\tau < b} \sigma_\T(\tau) \cdot (f \circ \sigma_\T(\tau) - f(\tau)) + \sum_{\tau < \sigma_\T(\tau)}^{\tau < b} [s \cdot f(s)]\Big|_{\tau}^{\sigma_\T(\tau)} \non \\
&&-\; \int_{\tau}^{\sigma_\T(\tau)} f \circ \sigma_{\Ttil}(s) \, \Delta_{\Ttil} s \non\\
& = & \int_{a}^{b} f \circ \sigma_{\Ttil}(s) \, \Delta_{\Ttil} s + \sum_{\tau < \sigma_\T(\tau)}^{\tau < b} \mu_\T(\tau) \cdot f(\tau) - \int_{\tau}^{\sigma_\T(\tau)} f \circ \sigma_{\Ttil}(s) \, \Delta_{\Ttil} s. \label{eq_id_final}
\end{eqnarray}


Note that since C$_{rd}^1([a,b]_\T,\R)$ is dense in C$_{rd}([a,b]_\T,\R)$ in the supremum norm, given $g \in$ C$_{rd}([a,b]_\T,\R)$ there is a sequence $f_n \in$ C$_{rd}^1([a, b]_\T,\R)$ such that, for each $f_n$ the identity (\ref{eq_id_final}) holds and $f_n$ converges to $g$ uniformly at almost every point at $[a,b] _{\Ttil}$. Therefore, the same identity holds for $g$. \hfill $\square$
\end{prova}


\bibliographystyle{unsrtnat}
\bibliography{references}  






\end{document}